# Graceful and Strongly Graceful Permutations

*Rafael Rofa*

**Abstract**
*A graceful labelling of a graph G* is an injective function f from the set of vertices of G into the set {0,1,…,|EG|} such that if edge uv is assigned the label |f(u)-f(v)| then all edge labels have distinct values. *A strong graceful labelling of a tree T with a perfect matching* is a graceful labelling of T with the additional property that the sum of the vertex labels of each odd labelled edge add up to |ET|. A lobster or a 2-distant tree is a tree T that contains a path P such that any vertex of T is a distance at most 2 from a vertex of P. In this paper, we define *generalised strongly graceful permutations* and discover two new permutations in addition to the known permutation that is obtained by replacing each vertex label f(v) by |ET| - f(v). We use these permutations to prove, by induction, that a lobster with a perfect matching that consists of the set of end edges of the lobster, is strongly graceful. Further, we show that there exist strongly graceful labellings that assign the label 0 to four specific vertices of any tree belonging to this family of lobsters. By using the technique developed in this paper we will, further, present a tractable way for proving an equivalent form of Bermond's conjecture which states that all lobsters are graceful. Two out of a total of three cases of the proposed equivalent form of Bermond's conjecture are completed leaving the third case open for refutation or completion.

1.  **Introduction and Background**

Throughout this paper, we only deal with simple graphs, that is graphs without loops and multiple edges and we use standard terminology of graph theory. A graph G is said to be a *graceful graph* if there is a injective function f: V(G) → {0, 1 . . ., |EG|} such that if edge uv is assigned the label |f(u) – f(v)| then all edge labels are pairwise distinct. Since no edge label can exceed |EG| this means that the edge labels are exactly the element of the set {1, 2, . . ., |EG|}. If v is in VG and f is a graceful labelling of G, then v and f(v) are used interchangeably to identify the vertex v. When applied to trees, which is the focus of this paper, the definition of graceful graphs is rephrased as follows: a tree is said to be graceful if there is a bijection f: VT → {0, 1, . . ., |ET|} such that if edge uv is assigned the label |f(u) – f(v)|, then the induced edge labels are exactly the elements of the set {1, 2, . . ., |ET|}.

The term graceful labelling was first coined by (Golomb, 1972) but the concept was first introduced by A. ROSA who instead, used the term β-valuations (Rosa, 1967). Rosa introduced β-valuations and three other types of graph labelling as tools for addressing the Ringel Conjecture that the complete graph $K_{2n+1}$ can be decomposed into 2n+1 isomorphic trees each on n edges (Ringel, 1963). The connection between the Ringel Conjecture and graph labelling lies in the fact that Ringel's conjecture is equivalent to the statement that all trees are graceful. This last statement, which was raised by Kotzig, Ringel, and Rosa is known as the graceful tree conjecture (GTC) and is one of the most challenging open problems in graph theory (Rosa, 1967).

There is an extensive list of papers on the subject of graceful trees. Efforts are directed towards proving that certain families of trees are graceful. For example, the following trees are known to be graceful: caterpillars (trees with the property that the removal of its end edges leaves a path); lobsters (trees with the property that the removal of its end edges leaves a caterpillar) with perfect matchings; symmetrical trees (rooted trees in which every level contains vertices of the same degree); trees with diameter at most 5; trees with at most four end vertices. Given that paths and caterpillars are graceful, one way to approach the GTC is to continue the progression and prove that lobsters, 3-distant trees and so on, are graceful. In 1979, Bermond conjectured that lobsters are graceful (Bermond, 1979). However even this conjecture is still open to this day. The reader is referred to the dynamic survey by J. Gallian (Gallian, 2018) for a detailed list of families of graceful graphs.

A tree T with a perfect matching M (and hence with an even number of vertices) is said to be *strongly graceful* if it has a graceful labelling f such that for every edge xy in M, f(x) + f(y) = |VT|-1. This is equivalent to saying that all odd edge labels are obtained from the absolute difference between the components of each pair of the set of vertex labels {{0, |VT|-1},{1,|VT|-2},{2,|VT|-3},…, {$\frac{|VT|}{2} - 1, \frac{|VT|}{2}$}}. Note that exactly one component of each pair is even. *The spike tree (spik(T))* of a tree on n vertices is obtained by adding n new vertices and n new edges to match every vertex of T. Therefore, by construction spik(T) is a tree that contains double the vertices of T and contains a perfect matching in which the set of edges that represent the perfect matching is exactly the set of end edges of spik(T). *The contree of a tree* with a perfect matching M is obtained by contracting every edge of M. Broersma and Hoede proved that the Graceful Tree Conjecture is equivalent to the conjecture that every tree containing a perfect matching is strongly graceful (Broersma & Hoede, 1999) The outline of the proof is as follows: Given a strongly graceful labelling of spik(T) for a tree T, a graceful labelling of T is obtained by giving labels to vertices of T according to half the even vertex label of the edge in M that has been contracted. Conversely, given a graceful labelling f of a tree T, a strongly graceful labelling of spik (T) is constructed by doubling every vertex label, giving the doubled label to the end vertices of M and assigning the matched vertices such that the sum adds to n-1. See figure 1 for an example (The bold edges are the edges of M). In particular, we note that because M is exactly the set of end edges of spik(T) then we have the following, more specific form of equivalence to the GTC.

**Theorem 1.** *A tree is graceful if and only if spik (T) is strongly graceful.*

Therefore, to prove the GTC, one does not have to prove that *every* tree with a perfect matching is strongly graceful, rather one has to show that every tree with a perfect matching M, where M is the set of end edges, is graceful.

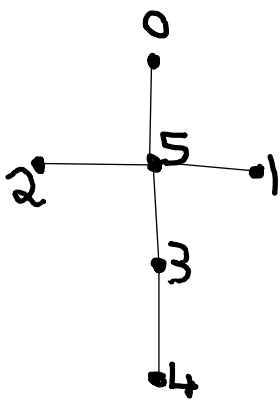 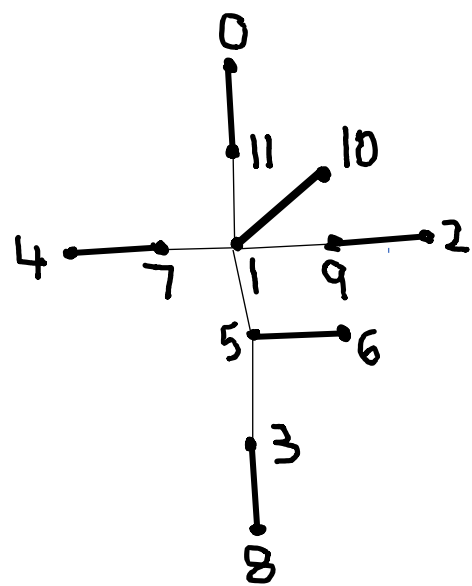

Figure 1(a). A graceful labelling of a tree T    Figure 1(b). A strong graceful labelling of Spik(T)

Therefore, we have the following equivalent form of the GTC

**Conjecture 1**. A tree T with a perfect matching that consists of all end edges of T is strongly graceful.

Since spik(L) of a lobster L is a 3-distant tree T with a perfect matching that consists of all end edges of T, then by Theorem 1, Bermond's conjecture that lobsters are graceful is equivalent to the following conjecture:

**Conjecture 2.** *A 3-distant tree T with a perfect matching M, where M is the set of all end edges of T is strongly graceful.*

Conjecture 1 can be generalised to:

**Conjecture 3**. *A (k+1)-distant tree T with a perfect matching M, where M is the set of all end edges of T and k ≤ Diameter(contree(T)) is strongly graceful.*

In this paper we prove Conjecture 3 for the case k = 1. The proof for the case k = 2 (conjecture 2), which is equivalent to Bermond's conjecture is attempted. However, only two out of three cases of the proof are completed leaving the third case open for refutation or completion. These proofs utilise the concept of generalised strongly graceful permutations that we introduced in the next section.

## 2. Graceful and strongly graceful permutations

**Definition 1.** Let T be a graceful (strongly graceful) tree on n vertices with graceful (strongly graceful) labelling f. A graceful (Strongly Graceful) permutation of VT with respect to f is a permutation g of the vertex labels 0,1, . . ., n-1 of f that preserves gracefulness (strong gracefulness). A generalised graceful (strongly graceful) permutation of T is a graceful (strongly graceful permutation) of T with respect to any graceful (strongly graceful) labelling. We denote the labelling that is obtained by applying the permutation g to the vertex labels of f by g[f]

These definitions could be extended to all graceful (strongly graceful) graphs, but we limit our definitions to trees in this paper.

**Example 1**. Given any graceful labelling (Strongly graceful labelling) f of T then
  a) the identity permutation e is a generalised graceful (strongly graceful) permutation of T.
  b) the permutation r = (0 n-1)(1 n-2)…($\frac{n}{2} - 1\ \frac{n}{2}$) if n is even and r = (0 n-1)(1 n-2)…($\frac{n-1}{2} - 1\ \frac{n-1}{2} + 1$)($\frac{n-1}{2}$) if n is odd ( g fixes the vertex label $\frac{n-1}{2}$) are generalised graceful ( generalised strongly graceful) permutations of T. Note that applying the permutation r to the vertex labels of f produces the known complimentary labelling of f that replaces each vertex label b by n-1-b so that the sum is n-1.

**Proof.** If uv∈VT, then |r(f(u)) - r(f(v))| = |n-1-f(u)-(n-1-f(v))| = |f(v) - f(u)| = |f(u) – f(v)|. Therefore, r preserves the edge labels of f edgewise and therefore, preserves gracefulness (strong gracefulness). Thus, r is a graceful (strongly graceful) permutation depending on whether f is a graceful or strongly graceful labelling of T. End of proof.

Note that a generalised graceful (strongly graceful) permutation is necessarily a graceful (strongly graceful) permutation, but the inverse is not necessarily true as the next example shows.

**Example 2.** g = (0 4 1 5)(2 3) is a graceful permutation of the tree T of figure 1(a) with respect to the graceful labelling, f that is shown in the figure. However, g is not a graceful permutation of T with respect the graceful labelling g[f]. Hence, g is not a generalised graceful labelling of T.

The only known generalised graceful (strongly graceful) permutations are the identity permutation and the permutation r of example 1(b) which is obtained by replacing each vertex label b by n-1-b. However, in theorems 2 and 3 below, we prove the existence of two additional generalised strong graceful permutations each consisting of a product of 2-cycles. Therefore, one would have to consider whether the set of generalised strong graceful permutations is closed under the operation of permutation composition. The following Lemma addresses this point, and it will be used in the proof of theorem 3.

**Lemma 1**. *If g and h are generalised strong graceful permutations of a tree T, then so is gh and hg (that is the set of generalised strong graceful permutations is closed under the operation of permutation product).*

Proof. Let f be a strong graceful labelling of a tree T, then g[f] and h[f] are both strong graceful labellings of T. Hence, gh[f] = g[h[f]] is a strong graceful labelling of T because h[f] is a strong graceful labelling of T and g is a generalised strong graceful labelling of T. Therefore, gh is a generalised graceful permutation of T. Because h and g are interchangeable then hg is also a generalised strong graceful permutation.

In the following theorems we find two additional generalised strong graceful permutations (other than the identity e and the permutation r of example 1) and use them in the proofs of Lemma 2 and conjecture 3 for k=1 and parts of the case k=2.

**Theorem 2.** *Let T be a strongly graceful tree on n vertices that contains a perfect matching defined by the set M of edges and let f be a strong graceful labelling of T, then the permutation $g_1$ = (0 1)(2 3) . . . (n-2 n-1) is a generalised strong graceful permutation of T.*

Note that $g_1$ increases every even labelled vertex, including 0, by 1 and decreases every odd vertex label by 1.

Proof. If uv ∉ M, then by definition, edge uv is assigned an even label under f. Since every even edge label must be obtained as the absolute difference between two vertex labels having the same parity (both even or both odd) and $g_1$ increases *every* even vertex label of f by 1 and decreases *every* odd vertex label by 1, then $g_1$ preserves all even edge labels edgewise.

On the other hand, if uv ∈ M, then by definition edge uv is assigned an odd label under f. Therefore, one of f(u) or f(v) is odd, and the other is even (including 0). Without loss of generality assume f(u) is even and f(v) is odd, then $g_1(f(u)) + g_1(f(v)) = f(u) + 1 + f(v) - 1 = f(u) + f(v) = n-1$. Therefore, $g_1$ preserves the condition of being strongly graceful and we have the following cases:

i) if $\frac{n}{2}$ is even, then $g_1$ interchanges the odd edge label between vertices f(u) and f(v) with the odd edge label between vertices f(u) +1 and f(v) -1.

ii) If $\frac{n}{2}$ is odd, then $g_1$ fixes odd edge label between vertices $\frac{n}{2} - 1$ and $\frac{n}{2}$ while interchanging all other odd labelled edges as in the case (i).

Therefore, $g_1$ is a strong graceful permutation as asserted by the theorem.

Note that given a strong graceful labelling f of T with f(u) =0 and f(v) = 1, f(w) =|V(T)|-2, and f(z) = |V(T)|-1 then, $g_1$ provides another strong graceful labelling of T, namely $g_1[f]$ with $g_1[f](v) = 0$, $g_1[f](u) = 1$, $g_1[f](w) = |V(T)|-1$, and $g_1[f](z) = |V(T)|-2$. That is the strong graceful labelling g[f] interchanges vertex labels 0 and 1 and also interchanges vertex labels |V(T)|-1 and |V(T)|-2 of f. This is going to be required in proving Lemma 2 and conjecture 3 for k=1 and parts of the case k=2.

The following example demonstrates $g_1$.

*Example 3*. Consider the strong graceful labelling f of the tree T, shown in Figure 2(a) below and consider the permutation $g_1$ = (0 1)(2 3)(4 5)(6 7)(8 9)(10 11), as defined in theorem 2. $g_1[f]$ is the strong graceful labelling of T that is shown in figure 2(b). Notice how $g_1$ [f] interchanges vertex labels 0,1 and vertex labels 10 and 11.

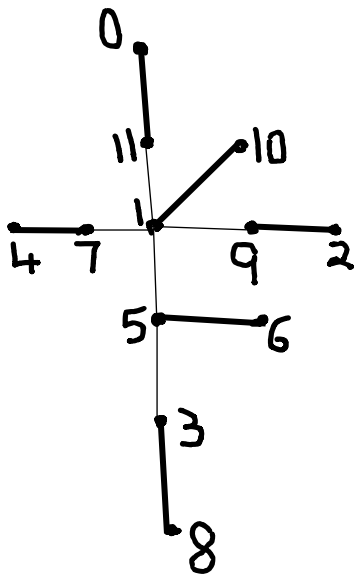
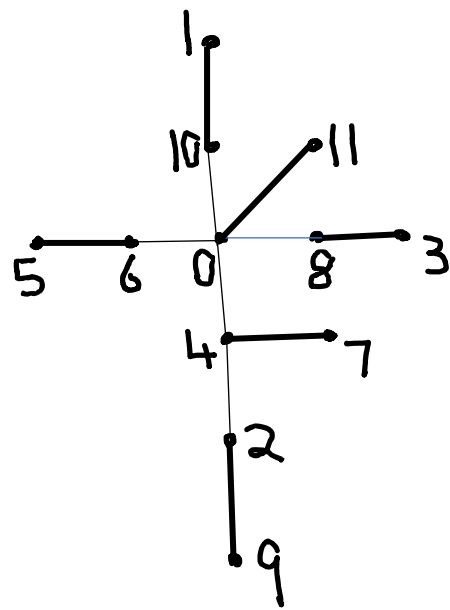

Figure 2(a)                                                        Figure 2(b)

**Theorem 3.** *Let T be a strongly graceful tree on n vertices that contains a perfect matching defined by the set M of edges and let f be a strong graceful labelling of T then*

i) *if $\frac{n}{2}$ is even, the* permutation $g_2 = (0\ n\text{-}2)(1\ n\text{-}1) \ldots (\frac{n}{2}-2\ \frac{n}{2})(\frac{n}{2}-1\ \frac{n}{2}+1)$ *is a generalised strong graceful permutation of T, and*

ii) *if $\frac{n}{2}$ is odd, the permutation $g_2 = (0\ n\text{-}2)(1\ n\text{-}1) \ldots (\frac{n}{2}-3\ \frac{n}{2}+1)(\frac{n}{2}-2\ \frac{n}{2}+2)(\frac{n}{2}-1)(\frac{n}{2})$ is a generalised strong graceful permutation of T.*

*Proof.* In both cases, $g_2$ is the product (in any order) of the generalised strong graceful permutations r (example 1) and $g_1$ (of Theorem 2). Hence, by Lemma 1, $g_2$ is a generalised strong graceful permutation of T as the theorem asserts.

$g_2$ can be thought of, and hence referred to as the complementary permutation to $g_1$.

Note that given a strong graceful labelling f of T with f(u) =0 and f(v) = 1, f(w) =|V(T)|-2, and f(z) = |V(T)|-1 then, $g_2$ provides another strong graceful labelling of T, namely $g_2[f]$ with $g_2[f](w) = 0$ and $g_2[f](u) = |V(T)|-2$, $g_2[f](z) = 1$, and $g_2[f](v) = |V(T)|-1$. That is the strong graceful labelling $g_2[f]$ interchanges vertex labels 0 and 1 and also interchanges vertex labels |V(T)|-1 and |V(T)|-2 of f. This is going to be required in proving Lemma 2 and conjecture 3 for k=1 and parts of the case k=2.

The following example demonstrates $g_2$.

**Example 4.** Consider the strong graceful labelling f of the tree T in Figure 2(a) above and consider the permutation $g_2 = (0\ 10)(1\ 11)(2\ 8)(3\ 9)(4\ 6)(5\ 7)$, as defined in Theorem 3. $g_2[f]$ is the strong graceful labelling of T that is shown in figure 3 below.

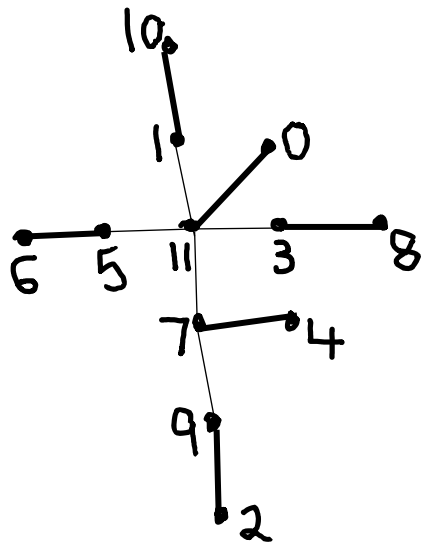

Figure 3

Notice that $\{e, r, g_1, g_2\}$ forms an abelian group under permutation composition where the inverse of each element is itself, and $rg_1= g_2$, $rg_2= g_1$ and $g_1g_2=r$. It is questionable whether there are other generalised strong graceful permutations that can be added to this group.

The following lemma utilises the generalised strongly graceful permutations r, $g_1$, and $g_2$ and it will be used in the proof of the case k=1 of conjecture 3 and part of the proof of the case k = 2 of the same conjecture.

**Lemma 2.** *Let T be a strongly graceful tree with a perfect matching M with $|VT|\geq 4$ and let f be a strong graceful labelling of T. Let $v_0$ be the vertex with $f(v_0) = 0$. Then*

(1) *there is a path $v_0, v_1, v_2, v_3$ in T with $v_0v_1 \in M$ and $v_2v_3 \in M$ such that $f(v_1) = |VT|-1$, $f(v_2) = 1$, $f(v_3) = |VT| - 2$, and*

(2) *there are strongly graceful labellings $f_1, f_2, f_3$ such that:*

    (i) $f_1(v_1) = 0, f_1(v_0) = |VT|-1, f_1(v_2) = |VT| -2, f_1(v_3) = 1$.

    (ii) $f_2(v_2) = 0, f_2(v_0) = 1, f_2(v_1) = |VT| -2, f_2(v_3) = |VT|-1$.

    (iii) $f_3(v_3) = 0, f_3(v_0) = |VT| -2, f_3(v_1) = 1, f_3(v_2) = |VT|-1$.

*Proof.* (1) Since T has a perfect matching M, then $v_0$ must be matched with a vertex say, $v_1$, so $v_0v_1 \in M$. Because $f(v_0) = 0$, and $v_0v_1 \in M$, then $f(v_1) = |VT|-1$. The only way to get edge label $|VT|-2$ is by having vertex $v_1$ joined to a vertex, say $v_2$ with $f(v_2) = 1$. Also, $v_2$ must be matched with a vertex, say $v_3$, so $v_2v_3 \in M$. Because $f(v_2) = 1$ and $v_2v_3 \in M$, it follows that $f(v_3) = |VT| - 2$. This proves part 1 of the lemma.

(2) The strongly graceful labellings $f_1, f_2$ and $f_3$ are obtained by applying the generalised strong graceful permutation r of example 1, $g_1$ of Theorem 2, and $g_2$ of Theorem 3, respectively, to f.

End of Proof.

### 3. Inductive proof of Conjecture 3 for the case k = 1

For the case k = 1, conjecture 3 asserts that a 2 distant tree T (that is a lobster) with a perfect matching M, where M is the set of all end edges of T is strongly graceful. This is proved by proving the following stronger result.

**Theorem 4.** *Let T be a 2-distant tree (lobster), on 4 or more (even) number of vertices with a perfect matching M where M is the set of all end edges of T. Let $P = v_0, v_1, v_2, v_3, \ldots v_p$ be a path of maximal length such that any vertex of T is at a distance of at most 2 from a vertex of P. Then there exist strongly graceful labellings f, $f_1, f_2$ and $f_3$ with $f(v_0) = 0, f_1(v_1) = 0, f_2(v_2) = 0$ and $f_3(u_2) = 0$, where $v_2u_2 \in M$.*

Before starting the inductive proof, we will list some properties of T including the identification of vertex $u_2$.

Lemma 3. If T is a lobster that satisfies the conditions of Theorem 4 then

i) $\deg(v_0) = 1$ and hence, $v_0v_1 \in M$.
ii) $\deg(v_1) = 2$
iii) Except for edges $v_0v_1$ and $v_{p-1}v_p$, each of the other edges of P do not belong to M.
iv) every vertex of P, other than $v_0, v_1, v_{p-1}$, and $v_p$, is matched with an end vertex that does not belong to VP. Let $u_i \notin VP$ be the vertex matched with vertex $v_i \in$ VP ($i \neq 0, 1, p-1,$ and $p$,

Proof

i) Follows from the fact that P is a maximal path.
ii) Follows from that $v_0v_1 \in M$ and P is a maximal path.
iii) Any such edge is not an end edge of T and hence does not belong to M.
iv) Follows from iii.

*Proof of Theorem 4 (by induction on |VT|):*

The only lobster T with $|VT| = 4$ that satisfies the conditions of theorem 4 is $P_3$ (= $v_0, v_1, v_2, v_3$). For this tree, f is the strongly graceful labelling with $f(v_0) = 0, f(v_1) = |VT|-1 = 3, f(v_2) = 1$, and $f(v_3) = |VT| -2 = 2$ and $f_1, f_2$ and $f_3$ are obtained by applying the generalised strong graceful permutations $r, g_1$ and $g_2$, respectively, to $f$. Therefore, the theorem is true for $|VT| = 4$. Next, we assume that the theorem is true for any lobster that satisfies the conditions of the theorem with $n-2$ or fewer (even) number of vertices. Now, let $T$ be a lobster that satisfies the condition of the Theorem with $|VT| = n$. To complete the induction proof, we have to show that T satisfies the theorem. For that purpose, we delete edges $v_0 v_1$ and $v_1 v_2$ from T to get a lobster T` with a perfect matching M`= M/{$v_0 v_1, v_1 v_2$} that satisfies the conditions of the theorem with $|VT`| = n-2$. Therefore, by the induction assumption T` satisfies the theorem. We have the following cases:

<u>case 1</u>. $\deg_{T`}(v_2) \geq 3$ and hence vertex $v_2$ is adjacent to the unique end vertex $u_2 \notin$ VP, with $v_2u_2 \in$ M` (as implied by point iv of lemma 3) and it is also joined to a number of 2 paths of the form $v_2xy$ not in P with edge $xy \in M$`. In this case, the path P` = $y,x, v_2,v_3 \ldots v_p$ is a path of maximal length in T`. Hence by the induction assumption there is a strong graceful labelling h of T` with $h(v_2) = 0$ and therefore $h(u_2) = n-3$ [because $v_2u_2 \in$ M`]. Now $v_2$ is joined to a number of 2 paths of the form $v_2xy$ not in P with edge $xy \in M$` *and because of the similarity of these path we can W.L.O.G assume that $h(x)=n-4$* [The only way to get the edge label n-4]. Also, since $xy \in M$` *and $h(x) = n-4$, then $h(y)=1$.* Now, reinstate edges $v_0 v_1$ and $v_1 v_2$ to obtain the original lobster T. Because $h(v_2) = 0$, then assigning the labels n-2 and -1 to vertices $v_1$ and $v_0$, respectively, then adding 1 to each vertex label (including labels of h) results in a graceful labelling f of T with $f(v_0) = 0, f(v_1) = n-1, f(v_2) = 1, f(u_2) = n-2$ and, $f(x) = n-3$, and $f(y) = 2$. To show that f is strongly graceful note that

(1) for edges $v_0v_1, v_2u_2,$ and $xy$ in M, *we have $f(v_0)+f(v_1)= f(v_2)+f(u_2)=f(x)+f(y)=n-1$.*

(2) for any other edge wz in *M*, we have f(w) + f(z) = h(w) + h(z) +2 = n-1 [because h is a strongly graceful labelling of T` with |VT`| = n-2 and f is obtained by adding 1 to each vertex label of h]. Therefore by (1) and (2) f is a strongly graceful labelling of T.

To complete the proof of this case we have to show that there are strong graceful labellings $f_1$, $f_2$ and $f_3$ of T such that $f_1(v_1) = 0$, $f_2(v_2) = 0$ and $f_3(u_2) = 0$. These are obtained by applying the generalised strong graceful permutations r, $g_1$ and $g_2$ from section 2, respectively, to f.

case 2. $Deg_{T`}(v_2) = 2$ and hence by part (iv) of Lemma 3, vertex $v_2$ is joined to a unique end vertex $u_2 \notin VP$ with $v_2u_2 \in M`$. In this case, the path P`= $u_2, v_2, v_3, v_4, \ldots v_p$ is a path of maximal length in T`. Hence by the induction assumption there is a strong graceful labelling h of T` with h($v_2$) = 0. Since $v_2u_2 \in M`$ and h is a strongly graceful labelling of T` then h($u_2$) = |VT`|-1 = n-3. Because $v_3$ is the only other vertex adjacent to $v_2$ (other than $u_2$), then h($v_3$) = 1. By lemma 1 part iv, $v_3$ is adjacent to a vertex $u_3 \notin VP$ with $v_3u_3 \in M`$, and therefore h($u_3$) = |VT`|-2 = n-4. Now, reinstate edges $v_0v_1$ and $v_1v_2$ to obtain lobster T. Because h($v_2$) = 0, then assigning the labels n-2 and -1 to vertices $v_1$ and $v_0$, respectively, then adding 1 to each vertex label (including labels of h) results in a graceful labelling f of T with *f($v_0$) = 0, f($v_1$) = n-1, f($v_2$) = 1,* f($v_3$) = 2, *f($u_2$) = n-2 and, f($u_3$) = n-3*. To show that f is strongly graceful note that

(1) for edges $v_0v_1, v_2u_2,$ and $v_3u_3$ in M, *we have f($v_0$) + f($v_1$) = f($v_2$) + f($u_2$) = f($v_3$) + f($u_3$) = n-1,* and

(2) for any other edge xy in *M*, we have f(x) + f(y) = h(x) + h(y) +2 = n-1 [h is a strongly graceful labelling of T` with |VT`| = n-2 and f is obtained by adding 1 to each vertex label of h]. Therefore by (1) and (2) f is a strongly graceful labelling of T.

To complete the proof of this case we have to show that there are strong graceful labellings $f_1$, $f_2$ and $f_3$ of T such that $f_1(v_1) = 0$, $f_2(v_2) = 0$ and $f_3(u_2) = 0$. These are obtained by applying the generalised strong graceful permutations r, $g_1$ and $g_2$ from section 2, respectively, to f.

End of proof.

## 4. Bermond's Conjecture

By using a technique that is similar to the one that was used in proving theorem 4, we now outline a proof for conjecture 2, which, as explained earlier, is equivalent to Bermond's conjecture. The inductive proof has 3 cases, and I was able to finalise 2 of them. I present the arguments here hoping that the third case would be completed or refuted in the future. The proof is attempted through the partial proof of the following stronger result.

Theorem 5. *Let T be a 3-distant tree on 4 or more (even) number of vertices with a perfect matching M where M is the set of all end edges of T. Let P = $v_0, v_1, v_2, v_3, \ldots v_p$ be a path of maximal length such that any vertex of T is at a distance of at most 3 from a vertex of P. Then there exist strongly graceful labellings f, $f_1$, $f_2$ and $f_3$ with f($v_0$) = 0, $f_1(v_1)$ = 0, $f_2(v_2)$ = 0 and $f_3(u_2) = 0$, where $v_2u_2 \in M$.*

Before starting the inductive proof, note that if T satisfies the conditions of theorem 5 then it must also satisfy statements i, ii, and iii of Lemma 3.

*Partial Proof of Theorem 5 (by induction on |VT|):*

The only 3-distant T with $|VT| = 4$ that satisfies the conditions of theorem 5 is $P_3$ ($= v_0, v_1, v_2, v_3$) which was shown in the proof of theorem 4 to satisfy the implication of Theorem 5. (Note that theorems 4 and 5 have the same implication which is the existence of a strong graceful labelling $f, f_1, f_2$ and $f_3$ with $f(v_0) = 0, f_1(v_1) = 0, f_2(v_2) = 0$ and $f_3(u_2) = 0$, where $v_2u_2 \in M$). Next, we assume that the theorem is true for any 3-distant tree that satisfies the conditions of the theorem with $n$-2 or fewer (even) number of vertices. Now, let $T$ be a 3-distant tree that satisfies the condition of the Theorem with $|VT| = n$. To complete the induction proof, we have to show that T satisfies the theorem. For that purpose, we delete edges $v_0 v_1$ and $v_1 v_2$ from T to get a 3-distant tree T` with a perfect matching $M` = M/\{v_0 v_1, v_1 v_2\}$ that satisfies the conditions of the theorem with $|VT`| = n$-2. Therefore, by the induction assumption T` satisfies the theorem. We have the following cases:

<u>case 1</u>. $Deg_{T`}(v_2) \geq 3$ and hence vertex $v_2$ is adjacent to the unique end vertex $u_2 \notin VP$, with $v_2u_2 \in M`$ (as implied by point iv of lemma 3) and it is also joined to a number of 2 paths of the form $v_2xy$ not in P with edge $xy \in M`$. In this case, the path P` $= y, x, v_2, v_3 \ldots v_p$ is a path of maximal length in T`. This is exactly the same as case 2 in the proof of theorem 4 and hence can be replicated here.

<u>case 2</u>. $Deg_{T`}(v_2) = 2$ and hence by part (iv) of Lemma 3, vertex $v_2$ is joined to a unique end vertex $u_2 \notin VP$ with $v_2u_2 \in M`$. We have the following subcases:

<u>case 2a</u>.     path P`$= u_2, v_2, v_3, v_4, \ldots v_p$ is a path of maximal length in T`. This is exactly the same as case 1 in the proof of theorem 4 and hence can be replicated here.

<u>Case 2b</u>.     path P` is not maximal. This can happen only if vertex $v_3$ is joined to a 3-path $v_3xyz$ whose edges are not in P and in this case the path $z, y, x, v_3, v_4, v_5, \ldots v_p$ is a path of maximal length. I am unable to apply the same techniques to prove that T satisfises the theorem in this case.

A possible way forward is to find additional generalised strong graceful permutations other than r, $g_1$, and $g_2$ that has been used in this paper or to prove that spik(L) of a lobster L is strongly graceful through proving a variant of theorem 5 that can be handled by the techniques that have been developed in this paper or other otherwise). If progress can be made on this case, then not only one would have completed the proof of Bermond's conjecture but would also have laid a foundation for proving that all 3-distant trees are graceful.